\documentclass[preprint,12pt]{elsarticle}
\usepackage{amssymb}
\usepackage{amsfonts}
\usepackage{indentfirst,latexsym,bm,amsmath,amsthm}
\usepackage{latexsym,mathrsfs,amsfonts}
\usepackage{bbm}
\usepackage{amsfonts}
\usepackage{indentfirst,latexsym,bm,amsmath,amsthm}
\usepackage{hyperref}
\makeatletter
\def\tank#1{\protected@xdef\@thanks{\@thanks
 \protect\footnotetext[0]{#1}}}
\def\bigfoot{

 \@footnotetext}
\makeatother

\newtheorem{remark}{Remark}[section]
\newtheorem{lemma}{Lemma}[section]

\newtheorem{theorem}{Theorem}[section]
\newtheorem{prop}{Proposition}[section]

\newtheorem{definition}{Defintion}[section]
\renewcommand{\theequation}{\arabic{section}.\arabic{equation}}

\begin{document}

\begin{frontmatter}
\title{\bf SPDEs with two reflecting walls and two singular drifts\thanks{This work was supported by the Fundamental Research Funds for the Central Universities, No.2013RC0906;
and by NSFC, No.11101419. }}
\author{Juan Yang$^a$,\ \ \   Jianliang Zhai$^b$}
\address[1]
{\  School of Science, Beijing University of posts and telecommunications, No.10 Xitucheng Road, Haidian District, Beijing, 100875, China. Email: yangjuanyj6@gmail.com}
\address[2]
{\  School of Mathematical Sciences, University of Science and Technology of China, Hefei, 230026, China. Email: zhaijianliang2008@163.com}

\begin{abstract}
We study SPDEs with two reflecting walls
$\Lambda^1$, $\Lambda^2$ and two singular drifts $\frac{c_1}{(X-\Lambda^1)^{\vartheta}}$, $\frac{c_2}{(\Lambda^2-X)^{\vartheta}}$, driven by space-time white noise. First, we establish
the existence and uniqueness of the solutions $X$ for $\vartheta\geq 0$. Second, we obtain the following pathwise properties of the solutions $X$.
If $\vartheta>3$, then a.s. $\Lambda^1<X<\Lambda^2$ for all $t\geq0$;
If $0<\vartheta<3$, then $X$ hits $\Lambda^1$ or $\Lambda^2$ with positive probability in finite time.
Thus $\vartheta=3$ is the critical parameter for $X$ to hit reflecting walls.



\end{abstract}

\end{frontmatter}

\section{Introduction} \setcounter{equation}{0}
\renewcommand{\theequation}{1.\arabic{equation}}

In this paper, we are concerned with the following reflected stochastic partial differential equations (SPDEs):
\begin{equation}\label{a01}
    \begin{cases}
\frac{\partial{X(x,t)}}{\partial{t}}=\frac{\partial^2{X(x,t)}}{\partial{x^2}}+f(x,t,X)+\frac{c_1}{(X-\Lambda^1)^{\vartheta}}- \frac{c_2}{(\Lambda^2-X)^{\vartheta}}
 \\
   \ \ \ \ \ \ \ \ \ \ \ +\chi\big(x,t,X\big)\dot{W}(x,t)+\Upsilon(x,t)-\Gamma(x,t),\ \ \ (x,t)\in{Q}:=[0,1]\times\mathbb{R}_+\\
    \Lambda^1(x,t)\leq{X}(x,t)\leq\Lambda^2(x,t).
                        \end{cases}
                        \end{equation}
Here ${W}$ is a space-time white noise and $\vartheta > 0$. Random measures $\Gamma$ and $\Upsilon$ are the additional forces
preventing $X$ from leaving $[\Lambda^1,\ \Lambda^2]$.

\vskip 0.3cm
Considering that SPDEs with reflection are natural extensions of the deterministic obstacle problems, there has been an upsurge of interest in
this topic.
The existence and uniqueness of the solutions have been well studied for the case of Lipchitz coefficient
with a single reflecting barrier
0, \emph{i.e.} $\Lambda^1=0,\Lambda^2=\infty, $ $c_1=c_2=0$, see Naulart and Pardoux \cite{np}, Donati-Martin and Pardoux \cite{dp1},  Xu and Zhang \cite{xz}, etc..
 Otobe \cite{os} and  Zhang and Yang \cite{zy} obtained the existence and uniqueness of the solutions of
SPDEs with double reflecting
barriers and Lipchitz coefficient, driven by an additive white noise and by multiplicative white noise respectively.

\vskip 0.3cm

All of the results mentioned above are devoted to the case
of Lipchitz coefficient. Mueller \cite{ml} proved that solutions of SPDEs $X_t=X_{xx}+g(X)\dot{W}(t,x)$, where  $g(X)$ is non Lipschiz, do not blow up in finite time, they did this by introducing a drift $X^{-\vartheta}$ into the former equation as a byproduct. It was proved that this drift forces solutions to stay positive with probability $1$.  Mueller and Pardoux \cite{mp} concentrated on the case when $0<\vartheta<3$ and showed that the solutions hit $0$ with positive probability.
Thus, $\vartheta=3$ is the critical exponent for $X$ to hit zero in finite time.
More precisely, if $\vartheta>3$, then $inf_x\ X(t,x)$ never reaches 0, and for $0<\vartheta<3$, $inf_x\ X(t,x)$
has a positive probability of reaching 0 in finite time. Existence and uniqueness for all time of a solution for $\vartheta=3$
was first verified by Zambotti \cite{zi}, and
it also showed the existence and uniqueness of the solution in the case of $0<\vartheta<3$ and a single reflecting barrier 0 in Zambotti \cite{zi}.
The hitting properties of the solution for $\vartheta=3$ was discussed in \cite{dmz}. All those discussions are under the assumption that SPDEs with the singular drift holds.

\vskip 0.3cm
Inspired by these, in this paper, we consider SPDEs
(\ref{a01}), which has double smooth reflecting walls $\Lambda^1$ and $\Lambda^2$ and two singular drifts $\frac{c_1}{(X-\Lambda^1)^{\vartheta}}$, $\frac{c_2}{(\Lambda^2-X)^{\vartheta}}$, for all $\vartheta>0$. As far as we know, there are few literatures that have studied the reflection problem with singularities.
To show the existence
and uniqueness of the solutions of Eq. \eqref{a01},
the strategy is
similar to that in Xu and Zhang \cite{xz} as well as \cite{zy}.
Due to the singular terms, we need to approximate the solution by introducing two monotone sequences.
As an extension from one wall to two walls case, we show that
for $\vartheta>3$, the solution of Eq. (\ref{a01}) stays in the interval
$\big(\Lambda^1,\ \Lambda^2\big)$ for all $t\geq 0$, a.s., that is $\Upsilon=\Gamma=0$;
for $0<\vartheta<3$, the solution of Eq. (\ref{a01}) hits $\Lambda^1$ or $\Lambda^2$ with positive probability in finite time.
Thus $\vartheta=3$ is the critical parameter for $X$ to hit reflecting walls. Various other
properties of the solutions have been studied, see
\cite{dmz}, \cite{dp2}, \cite{yz}, \cite{za} and \cite{zw}.
\vskip 0.3cm

The paper is organized as follows.
Definitions and assumptions are given in Section 2. Section 3 is to establish
the existence and uniqueness of the solution of Eq. \eqref{a01}. In Section 4 we
consider the pathwise properties of Eq. \eqref{a01},
Subsection 4.1 is devoted to the case $\vartheta>3$, and in
Subsection 4.2, we deal with the case when $0<\vartheta<3$.

\section{Preliminaries}
\setcounter{equation}{0}
\renewcommand{\theequation}{2.\arabic{equation}}

In this section, we state the precise assumptions on the coefficients and the concept of solution.

First, we formulate Eq. \eqref{a01} as follows:
\begin{equation}\label{a02}
\frac{\partial{X}}{\partial{t}}=\frac{\partial^2{X}}{\partial{x^2}}+f(x,t,X)+\frac{c_1}{(X-\Lambda^1)^{\vartheta}}- \frac{c_2}{(\Lambda^2-X)^{\vartheta}}+\chi(x,t,X)\dot{W}+\Upsilon-\Gamma;
\end{equation}
with the conditions
\begin{equation}\label{a02}
\left \{\begin{array}{ll}
        X(0,t)=0,\ X(1,t)=0,\ \ {\rm for}\ t\geq0;\\
        X(x,0)=X_0(x)\in{C}([0,1]);\\
        \Lambda^1(x,t)\leq{X}(x,t)\leq{\Lambda^2}(x,t),\ \ {\rm for}\ (x,t)\in{Q}:=[0,1]\times\mathbb{R}_+.
\end{array}
\right.
\end{equation}

Let $f, \chi:[0,1]\times\mathbb{R}_+\times{C}(Q)\rightarrow\mathbb{R}$ be given measurable maps. Set $\mathcal{A}_{M,T}=\{X\in C(Q):\ \underset{x\in[0,1],t\in[0,T]}{\sup}|X(x,t)|\leq{M}\}$.
Introduce the following conditions:\\
(H0) $\Lambda^i(x,t)\in C(Q)$ satisfy
$\Lambda^i(1,t)\leq{0}$,  $\Lambda^i(1,t)\geq{0}$, for $i=1,2$;\\
(H1)  $\Lambda^1(x,t)<\Lambda^2(x,t)$, $\forall(x,t)\in(0,1)\times \mathbb{R}_+$; \\
(H2)
$\frac{\partial{\Lambda^i}}{\partial{t}}+\frac{\partial^2{\Lambda^i}}{\partial{x^2}}\in
L^2([0,1]\times[0,T])$;\\
(H3)
$\frac{\partial{}}{\partial{t}}\Lambda^i(0,t)=\frac{\partial{}}{\partial{t}}\Lambda^i(1,t)=0$
for $t\geq0$;\\
(H4) $\frac{\partial{}}{\partial{t}}(\Lambda^2-\Lambda^1)\geq0$.\\
(F1) $f(\cdot,\cdot;0),\chi(\cdot,\cdot;0)\in L^2([0,1]\times[0,T])$, for $T>0$;\\
(F2) for every $X,\hat{X}\in\mathcal{A}_{M,T}$ there exists $C_{T,M}>0$ such that
$$|f(x,t;X)-f(x,t;\hat{X})|+|\chi(x,t;X)-\chi(x,t;\hat{X})|\leq{C_{T,M}}\underset{y\in[0,1],s\in[0,t]}{\sup}|{X}(y,s)-\hat{X}(y,s)|$$
for all $(x,t)\in[0,1]\times[0,T]$;\\
(F3) there exists $C_T>0$ such that
$$|f(x,t;X)|+|\chi(x,t;X)|\leq{C_T}(1+\underset{y\in[0,1],s\in[0,t]}{\sup}|X(y,s)|),\ \forall (x,t)\in[0,1]\times[0,T].$$

Now we recall the concept of solution of Eq. (\ref{a02}) in \cite{yz}.
\begin{definition} A triplet
$(X,\Upsilon,\Gamma)$ defined on a filtered probability space\\
$(\Omega,P,\mathcal {F};\{\mathcal {F}_t\})$ is a solution to the
SPDE \eqref{a02}, if\\
(i) $X\in{C(Q)}$ is a adapted random field
satisfying $\Lambda^1(x,t)\leq{X}(x,t)\leq{\Lambda^2}(x,t)$,
$X(0,t)=0$, and $X(1,t)=0$, a.s;\\
(ii) $\Upsilon(dx,dt)$ and $\Gamma(dx,dt)$ are positive and adapted random measures on $(0,1)\times\mathbb{R}_+$
satisfying
$$\Upsilon\big((\epsilon,1-\epsilon)\times[0,T]\big)<\infty,\ \Gamma\big((\epsilon,1-\epsilon)\times[0,T]\big)<\infty$$
for $\epsilon\in(0,\frac{1}{2})$ and $T>0$;\\
(iii) for every $t\geq0$ and $\psi\in{C_0^\infty}(0,1)$(the set of smooth
functions with compact supports)
\begin{eqnarray}\label{a03}
&&\big(X(t),\psi\big)-\big(X_0,\psi\big)-\int_0^t(X(s),\psi^{''})ds-\int_0^t\big(f(y,s,X),\psi\big)ds \nonumber\\
&=&\int_0^t\big(\frac{c_1}{(X-\Lambda^1)^{\vartheta}}
-\frac{c_2}{(\Lambda^2-X)^{\vartheta}},\psi\big)ds+\int_0^t\int_0^1(\chi(y,s,X),\psi)  W(dx,ds) \nonumber\\
&&+\int_0^t\int_0^1\psi\Upsilon(dx,ds)-\int_0^t\int_0^1\psi\Gamma(dx,ds),\ a.s,
\end{eqnarray}
where $(,)$ denotes the inner product in $L^2([0,1])$ and $X(t)$ denotes $X(\cdot,t)$;\\
(iv)$$\int_Q\big(X(x,t)-\Lambda^1(x,t)\big)\Upsilon(dx,dt)=\int_Q\big(\Lambda^2(x,t)-X(x,t)\big)\Gamma(dx,dt)=0.$$
\end{definition}

\section{The existence and uniqueness of solutions}
\setcounter{equation}{0}
\renewcommand{\theequation}{3.\arabic{equation}}
In this section, we prove the following theorem.

\begin{theorem}\label{t1}  Let $X_0\in{C}([0,1])$ satisfy
$\Lambda^1(x,0)\leq{X_0}(x)\leq{\Lambda^2}(x,0),\ X_0(0)=X_0(1)=0$. Assume conditions (H0)-(H4), (F1)-(F3) hold.
Then there exists a unique solution $(X, \Upsilon,\Gamma)$ to Eq. \eqref{a02}
 for all $\vartheta\geq0$.
\end{theorem}

As a first step we show that, for every $v(x,t)$ (sometimes denote by $v$) $\in{}C(Q)$, $v(x,0)=X_0(x)$ and
$v(0,t)=0,\ v(1,t)=0$, the deterministic PDE
\begin{equation}\label{a0001}
\left \{\begin{array}{ll}
       \frac{\partial{\Xi(x,t)}}{\partial{t}}=\frac{\partial^2{\Xi(x,t)}}{\partial{x^2}}+\frac{c_1}{(\Xi+v-\Lambda^1)^{\vartheta}}- \frac{c_2}{(\Lambda^2-\Xi-v)^{\vartheta}}+\Upsilon(x,t)-\Gamma(x,t);\\
        \Xi(0,t)=\Xi(1,t)=0,\ \ {\rm for}\ t\geq0;\\
        \Xi(x,0)=0,\ \ {\rm for}\ x\in[0,1];\\
        \Lambda^1(x,t)\leq{\Xi}(x,t)+v(x,t)\leq{\Lambda^2}(x,t),\ \ {\rm for}\ (x,t)\in{Q}
\end{array}
\right.
\end{equation}
has a unique continuous solution.

We now recall the precise definition of the solution for Eq. (\ref{a0001}) in \cite{yz}.
\begin{definition} A triplet $(\Xi,\Upsilon,\Gamma)$ is called a
solution to the PDE (\ref{a0001}) if\\
(i) $\Xi={\Xi(x,t);(x,t)\in{Q}}$ is a continuous function satisfying
$\Lambda^1(x,t)\leq{\Xi}(x,t)+v(x,t)\leq{\Lambda^2}(x,t),\ \Xi(x,0)=0,\
\Xi(0,t)=\Xi(1,t)=0$;\\
(ii) $\Upsilon(dx,dt)$ and $\Gamma(dx,dt)$ are measures on
$(0,1)\times\mathbb{R}_+$ satisfying
$$\Upsilon\big((\epsilon,1-\epsilon)\times[0,T]\big)<\infty,\  \Gamma\big((\epsilon,1-\epsilon)\times[0,T]\big)<\infty$$
for every $\epsilon\in(0,\frac{1}{2})$ and $T>0$;\\
(iii) for all $t\geq0$ and $\psi\in{C_0^\infty}(0,1)$ we have
\begin{eqnarray}\label{3.2}
     &&\big(\Xi(t),\psi\big)-\int_0^t(\Xi(s),\psi^{''})ds-\int_0^t\big(\frac{c_1}{(\Xi+v-\Lambda^1)^{\vartheta}}
-\frac{c_2}{(\Lambda^2-\Xi-v)^{\vartheta}},\psi\big)ds \nonumber\\
     &=&\int_0^t\int_0^1\psi\Upsilon(dx,ds)-\int_0^t\int_0^1\psi\Gamma(dx,ds),
\end{eqnarray}
where $\Xi(t)$ denotes $\Xi(\cdot,t)$.\\
(iv)
\begin{eqnarray*}
     &&\int_Q\big(\Xi(x,t)+v(x,t)-\Lambda^1(x,t)\big)\Upsilon(dx,dt) \\
          &=&\int_Q\big(\Lambda^2(x,t)-\Xi(x,t)-v(x,t)\big)\Gamma(dx,dt)\\
          &=&0.
\end{eqnarray*}
\end{definition}

\begin{prop}\label{p1}  Assume conditions (H0)-(H4) hold. Then
there exists a unique solution $(\Xi, \Upsilon,\Gamma)$ to Eq.
\eqref{a0001}.
Moreover,
denote $(\hat{\Xi},\hat\Upsilon,\hat\Gamma)$ by the solution to Eq. (\ref{a0001}) replacing $v$
with $\hat{v}$. Then
\begin{eqnarray}\label{a0002}
\|\Xi-\hat{\Xi}\|_\infty^T\leq\|v-\hat{v}\|_\infty^T,\ \forall T>0,
\end{eqnarray}
where
$\|J\|_\infty^T:=\underset{x\in[0,1],t\in[0,T]}{\sup}|J(x,t)|$.
\end{prop}

Before the proof of Proposition \ref{p1}, the following lemma for comparison is needed.
\begin{lemma}\label{l1}  Assume (H0)-(H4). Let $g^{\rho }:[0,1]\times \mathbb{R}\rightarrow \mathbb{R}$ be measurable,
$\rho >0$, such that $y\mapsto g^{\rho}(x,y)$ is monotone nonincreasing, Lipschitz-continuous uniformly in $x \in [0,1]$
and satisfies
$$|g^{\rho}(x,y)| \leq C(1+|y|),\ \ {\rm for}\  y\in \mathbb R, \ \rho>0,$$
for some $C\geq 0$. Consider the following PDE with reflection:
\begin{equation}\label{a05}
\left \{\begin{array}{ll}
        \frac{\partial{\Xi^\rho}}{\partial{t}}=\frac{\partial^2{\Xi^\rho}}{\partial{x^2}}+g^\rho(\Xi^{\rho})+\Upsilon^{\rho}-\Gamma^{\rho};\\
        \Xi^\rho(0,t)=\Xi^\rho(1,t)=0,\ \ {\rm for}\ t\geq0;\\
        \Xi^\rho(x,0)=0,\ \ {\rm for}\ x\in[0,1],
\end{array}
\right.
\end{equation}
where $\Upsilon^{\rho}$, $\Gamma^{\rho}$ satisfy $$\int_Q\big(\Xi^{\rho}(x,t)+v(x,t)-\Lambda^1(x,t)\big)\Upsilon^{\rho}(dx,dt)=0,$$
$$\int_Q\big(\Lambda^2(x,t)-\Xi^{\rho}(x,t)-v(x,t) \big)\Gamma^{\rho}(dx,dt)=0.$$
Then there exists a unique solution of Eq. \eqref{a05}, denoted by $\Xi^{\rho}$.
Moreover, if $\rho\mapsto g^{\rho}$ is monotone nonincreasing, then
$\rho\mapsto \Xi^{\rho}$ is monotone nonincreasing.
\end{lemma}

\begin{remark}\label{r1} To get the monotonicity of $\Xi^{\rho}$ w.r.t. $\rho$, the monotonicity of drift $g^{\rho}$ in
Eq. \eqref{a05} w.r.t. $\rho$ is required. If we assume that there exists the solution of Eq. \eqref{a05},
the Lipschitz-continuous and linear growth conditions of $g$ are not needed in Lemma \ref{l1}.
\end{remark}
\noindent\textbf{Proof.}\quad
The existence of the solution of Eq. \eqref{a05} can be obtained from the convergence of the penalized equation which is similarly to Theorem 2.1 in Otobe \cite{os}. The uniqueness of the solutions of \eqref{a05} can be proved in a similar way as Theorem 3.1 in Zhang and Yang \cite{zy}.
\vskip 0.3cm

Now we prove the comparison principle. Let $\rho_1\geq \rho_2>0$ and set $\phi=(\Xi^{\rho_1}- \Xi^{\rho_2})^+$, where $X^{\rho_1},\ X^{\rho_2}$ are the solutions of Eq. \eqref{a05} with drifts $g^{\rho_1}$,
$g^{\rho_2}$ respectively.
Consider the inner product between $\phi$ and $\Xi^{\rho_1}- \Xi^{\rho_2}$ in $L^2(0,1)$.
We have
\begin{eqnarray*}
     &&\| \phi(t) \|^2_{L^2}\\
     &=&2 \int_0^t <\phi,\ \frac{\partial^2}{\partial x^2}( \Xi^{\rho_1}-  \Xi^{\rho_2}  )>ds  +2 \int_0^t <\phi,\ g^{\rho_1}(\Xi^{\rho_1})-  g^{\rho_2}(\Xi^{\rho_2})>ds \nonumber\\
     && +  2 \int_0^t  <\phi,\ d\Upsilon^{\rho_1}- d\Upsilon^{\rho_2}  >  - 2 \int_0^t  <\phi,\ d\Gamma^{\rho_1}- d\Gamma^{\rho_2}  >  \nonumber\\
      &= & - 2 \int_0^t\| \frac{\partial \phi}{\partial x} \|^2_{L^2}ds + 2 \int_0^t <\phi,\ g^{\rho_1}(\Xi^{\rho_1})-  g^{\rho_1}(\Xi^{\rho_2})>ds\\
    && + 2 \int_0^t <\phi,\ g^{\rho_1}(\Xi^{\rho_2})-  g^{\rho_2}(\Xi^{\rho_2})>ds   \nonumber\\
     && +  2 \int_0^t  {1}_{\{\Xi^{\rho_1}\geq \Xi^{\rho_2}\}  } <\Xi^{\rho_1}+v-\Lambda^1- (\Xi^{\rho_2}+v-\Lambda^1)  ,\ d\Upsilon^{\rho_1}
    - d\Upsilon^{\rho_2}  > \\
     &&- 2 \int_0^t{1}_{\{\Xi^{\rho_1}\geq \Xi^{\rho_2} \} } <\Xi^{\rho_1}+v-\Lambda^2- (\Xi^{\rho_2}+v-\Lambda^2) ,\ d\Gamma^{\rho_1}- d\Gamma^{\rho_2}  > \\
      &=&  - 2 \int_0^t\| \frac{\partial \phi}{\partial x} \|^2_{L^2}ds + 2 \int_0^t <\phi,\ g^{\rho_1}(\Xi^{\rho_1})-  g^{\rho_1}(\Xi^{\rho_2})>ds\\
     &&+ 2 \int_0^t <\phi,\ g^{\rho_1}(\Xi^{\rho_2})-  g^{\rho_2}(\Xi^{\rho_2})>ds   \nonumber\\
     && +  2 \int_0^t {1}_{\{\Xi^{\rho_1}\geq \Xi^{\rho_2}\}  } (<\Xi^{\rho_1}+v-\Lambda^1  ,\ - d\Upsilon^{\rho_2}  > +<- (\Xi^{\rho_2}+v-\Lambda^1)  ,\ d\Upsilon^{\rho_1}  > )\\
      &&- 2 \int_0^t {1}_{\{  \Xi^{\rho_1}\geq \Xi^{\rho_2} \} }( <\Xi^{\rho_1}+v-\Lambda^2 ,\ - d\Gamma^{\rho_2}  >
      +<- (\Xi^{\rho_2}+v-\Lambda^2) ,\ d\Gamma^{\rho_1}  > )  \nonumber\\
     &\leq&  2 \int_0^t <\phi,\ g^{\rho_1}(\Xi^{\rho_2})-  g^{\rho_2}(\Xi^{\rho_2})>ds.
\end{eqnarray*}
\hfill$\Box$

\noindent\textbf{Proof of Proposition \ref{p1}}\quad  Existence: Consider the following reflected PDE with Lipschitz coefficients:
\begin{equation}\label{a04}
\left \{\begin{array}{ll}
       \frac{\partial{\Xi^{\epsilon_1,\epsilon_2}}}{\partial{t}}=\frac{\partial^2{\Xi^{\epsilon_1,\epsilon_2}}}{\partial{x^2}}
       +\frac{c_1}{[\epsilon_1+ (\Xi^{\epsilon_1,\epsilon_2}+v-\Lambda^1) ]^{\vartheta}}- \frac{c_2}{[\epsilon_2+ (\Lambda^2-\Xi^{\epsilon_1,\epsilon_2}-v) ]^{\vartheta}} +\Upsilon^{\epsilon_1,\epsilon_2}-\Gamma^{\epsilon_1,\epsilon_2};\\
       \Xi^{\epsilon_1,\epsilon_2}(0,t)=\Xi^{\epsilon_1,\epsilon_2}(1,t)=0;\\
       \Xi^{\epsilon_1,\epsilon_2}(x,0)=0,
\end{array}
\right.
\end{equation}
where $\epsilon_1,\epsilon_2>0$ and $\Upsilon^{\epsilon_1,\epsilon_2}$, $\Gamma^{\epsilon_1,\epsilon_2}$ satisfy $\int_Q\big(\Xi^{\epsilon_1,\epsilon_2}(x,t)+v(x,t)-\Lambda^1(x,t)\big) \\  \Upsilon^{\epsilon_1,\epsilon_2}(dx,dt)=0$, $\int_Q\big(\Lambda^2(x,t)-\Xi^{\epsilon_1,\epsilon_2}(x,t)-v(x,t)\big)\Gamma^{\epsilon_1,\epsilon_2}(dx,dt)=0$.
\vskip 0.3cm

The existence and uniqueness of the solutions of Eq. \eqref{a04} can be got from Lemma \ref{l1}.
Indeed, there is a unique solution ${\tilde{\Xi}}^{\epsilon_1,\epsilon_2}\in [\Lambda^1,\ \Lambda^2]$ of the following PDE
\begin{equation}\label{a041}
\left \{\begin{array}{ll}
       \frac{\partial{{\tilde{\Xi}}^{\epsilon_1,\epsilon_2}}}{\partial{t}}=\frac{\partial^2{{\tilde{\Xi}}^{\epsilon_1,\epsilon_2}}}{\partial{x^2}}
       +\frac{c_1}{[\epsilon_1+ ({\tilde{\Xi}}^{\epsilon_1,\epsilon_2}+v-\Lambda^1) \vee 0 ]^{\vartheta}}- \frac{c_2}{[\epsilon_2+ (\Lambda^2-{\tilde{\Xi}}^{\epsilon_1,\epsilon_2}-v)  \vee 0 ]^{\vartheta}} +\tilde\Upsilon^{\epsilon_1,\epsilon_2}-\tilde\Gamma^{\epsilon_1,\epsilon_2};\\
       {\tilde{\Xi}}^{\epsilon_1,\epsilon_2}(0,t)={\tilde{\Xi}}^{\epsilon_1,\epsilon_2}(1,t)=0;\\
       {\tilde{\Xi}}^{\epsilon_1,\epsilon_2}(x,0)=0,
\end{array}
\right.
\end{equation}
where $\tilde\Upsilon^{\epsilon_1,\epsilon_2}$, $\tilde\Gamma^{\epsilon_1,\epsilon_2}$ satisfy $\int_Q\big(\tilde \Xi^{\epsilon_1,\epsilon_2}(x,t)+v(x,t)-\Lambda^1(x,t)\big)  \tilde\Upsilon^{\epsilon_1,\epsilon_2}(dx,dt)=0$, $\int_Q\big(\Lambda^2(x,t)-\tilde \Xi^{\epsilon_1,\epsilon_2}(x,t)-v(x,t)\big)\tilde\Gamma^{\epsilon_1,\epsilon_2}(dx,dt)=0$.
Therefore, Eq.s \eqref{a04} and \eqref{a041} have the same solution.
Denote the solution of Eq. \eqref{a04} by $\Xi^{\epsilon_1,\epsilon_2}$.
\vskip 0.3cm

Fix $\epsilon_2$, by Lemma \ref{l1} we get that $\epsilon_1\mapsto \Xi^{\epsilon_1,\epsilon_2}$ is nonincreasing. Notice that
$\omega^{\epsilon_1,\epsilon_2}:=\epsilon_1 + \Xi^{\epsilon_1,\epsilon_2}$ is solution of
\begin{equation}\label{a07}
\left \{\begin{array}{ll}
       \frac{\partial{\omega^{\epsilon_1,\epsilon_2}}}{\partial{t}}=\frac{\partial^2{\omega^{\epsilon_1,\epsilon_2}}}{\partial{x^2}}
       +\frac{c_1}{ (\omega^{\epsilon_1,\epsilon_2}+v-\Lambda^1)^{\vartheta}}- \frac{c_2}{(\Lambda^2-\omega^{\epsilon_1,\epsilon_2}-v+  \epsilon_1+\epsilon_2)^{\vartheta}} +\Upsilon_\omega^{\epsilon_1,\epsilon_2}-\Gamma_\omega^{\epsilon_1,\epsilon_2};\\
       \omega^{\epsilon_1,\epsilon_2}(0,t)=\omega^{\epsilon_1,\epsilon_2}(1,t)=\epsilon_1;\\
       \omega^{\epsilon_1,\epsilon_2}(x,0)=\epsilon_1,
\end{array}
\right.
\end{equation}
where $\Upsilon_\omega^{\epsilon_1,\epsilon_2}$, $\Gamma_\omega^{\epsilon_1,\epsilon_2}$ satisfy $\int_Q\big(\omega^{\epsilon_1,\epsilon_2}(x,t)+v(x,t)-\Lambda^1(x,t)\big)\Upsilon_\omega^{\epsilon_1,\epsilon_2}(dx,dt)=0$, $\int_Q\big(\Lambda^2(x,t)-\omega^{\epsilon_1,\epsilon_2}(x,t)-v(x,t)\big)\Gamma_\omega^{\epsilon_1,\epsilon_2}(dx,dt)=0$.
\vskip 0.3cm

Let $g^{\epsilon_1}=\frac{c_1}{ (\omega^{\epsilon_1,\epsilon_2}+v-\Lambda^1)^{\vartheta}}- \frac{c_2}{(\Lambda^2-\omega^{\epsilon_1,\epsilon_2}-v+  \epsilon_1+\epsilon_2)^{\vartheta}}$. By Lemma \ref{l1}, we have for all $\epsilon_1''\geq \epsilon_1'>0$,
\begin{eqnarray}
     && \| (\omega^{\epsilon_1',\epsilon_2}-\omega^{\epsilon_1'',\epsilon_2}  )^+  \|^2_{L^2} \nonumber\\
     &\leq & \int_0^t\int_0^1[g^{\epsilon_1'}(  \omega^{\epsilon_1'',\epsilon_2}  )   -  g^{\epsilon_1''}(  \omega^{\epsilon_1'',\epsilon_2}  )](\omega^{\epsilon_1',\epsilon_2}-\omega^{\epsilon_1'',\epsilon_2}  )^+  dxds  \nonumber\\
     &= &\int_0^t\int_0^1 [  \frac{c_1}{ (\omega^{\epsilon_1'',\epsilon_2}+v-\Lambda^1)^{\vartheta}}- \frac{c_2}{(\Lambda^2-\omega^{\epsilon_1'',\epsilon_2}-v+  \epsilon_1'+\epsilon_2)^{\vartheta}} \nonumber\\
       &&- ( \frac{c_1}{ (\omega^{\epsilon_1'',\epsilon_2}+v-\Lambda^1)^{\vartheta}}- \frac{c_2}{(\Lambda^2-\omega^{\epsilon_1'',\epsilon_2}-v+  \epsilon_1''+\epsilon_2)^{\vartheta}}  )     ] \nonumber\\
       && (\omega^{\epsilon_1',\epsilon_2}-\omega^{\epsilon_1'',\epsilon_2}  )^+  dxds      \nonumber\\
        &\leq &0.
\end{eqnarray}

This means $\epsilon_1\mapsto\epsilon_1+ \Xi^{\epsilon_1,\epsilon_2}$ is nondecreasing. By the definition of \emph{uniformly continuous}, we obtain that $\Xi^{\epsilon_1,\epsilon_2}$ converges uniformly on $[0,1]\times[0, T]$ to a continuous function, denoted by $\Xi^{\epsilon_2}$, as $\epsilon_1 \downarrow 0$.
Moreover, $\epsilon_1\mapsto  \Upsilon^{\epsilon_1,\epsilon_2}-\Gamma^{\epsilon_1,\epsilon_2}$ is monotone nondecreasing and then $\Upsilon^{\epsilon_1,\epsilon_2}-\Gamma^{\epsilon_1,\epsilon_2}$ converges distributionally to some measure $\Upsilon^{\epsilon_2}-\Gamma^{\epsilon_2}$. From
the proof of existence of the solution in Theorem 3.1 in \cite{zy}, $\Upsilon^{\epsilon_1,\epsilon_2}$, $\Gamma^{\epsilon_1,\epsilon_2}$ converge distributionally to two measures $\Upsilon^{\epsilon_2}$, $\Gamma^{\epsilon_2}$ respectively, and denoted by $(\Xi^{\epsilon_2},\Upsilon^{\epsilon_2},\Gamma^{\epsilon_2} ) $ the solution of the following equation:
\begin{equation}\label{a08}
\left \{\begin{array}{ll}
       \frac{\partial{\Xi^{\epsilon_2}}}{\partial{t}}=\frac{\partial^2{\Xi^{\epsilon_2}}}{\partial{x^2}}
       +\frac{c_1}{ (\Xi^{\epsilon_2}+v-\Lambda^1)^{\vartheta}}- \frac{c_2}{[\epsilon_2+ (\Lambda^2-\Xi^{\epsilon_2}-v) ]^{\vartheta}}+\Upsilon^{\epsilon_2}-\Gamma^{\epsilon_2};\\
       \Xi^{\epsilon_2}(0,t)=\Xi^{\epsilon_2}(1,t)=0;\\
       \Xi^{\epsilon_2}(x,0)=0
\end{array}
\right.
\end{equation}
where $\Upsilon^{\epsilon_2}$, $\Gamma^{\epsilon_2}$ satisfy $$\int_Q\big(\Xi^{\epsilon_2}(x,t)+v(x,t)-\Lambda^1(x,t)\big)\Upsilon^{\epsilon_2}(dx,dt)=0,$$ $$\int_Q\big(\Lambda^2(x,t)-\Xi^{\epsilon_2}(x,t)-v(x,t)\big)\Gamma^{\epsilon_2}(dx,dt)=0.$$

In view of Remark 
3.1
, we could again apply Lemma \ref{l1} to Eq. \eqref{a08}. It is then easy to get that
 $\epsilon_2\mapsto \Xi^{\epsilon_2}$ is nonincreasing. $\epsilon_2\mapsto\epsilon_2+ \Xi^{\epsilon_2}$ is also nondecreasing,
and thus $\Xi^{\epsilon_2}$ converges uniformly as $\epsilon_2 \downarrow 0$ on $[0,1]\times[0, T]$ to a continuous function denoted by $\Xi$.
Also, $\Upsilon^{\epsilon_2}$, $\Gamma^{\epsilon_2}$ converge distributionally to two measures $\Upsilon$, $\Gamma$, respectively, and $(\Xi,\Upsilon,\Gamma ) $ is the solution of Eq.
\eqref{a0001}.
\vskip 0.3cm

Uniqueness: suppose $(\Xi_1,\Upsilon_1,\Gamma_1)$ and $(\Xi_2,\Upsilon_2,\Gamma_2)$
are two solutions to Eq. \eqref{a0001}. 
We can
choose an appropriate test function $\phi$ in Eq. (\ref{3.2}), following the same arguments as that in section 2.3 of Nualart and Pardoux \cite{np} and consequently get $\Xi_1=\Xi_2$. Using the same method of separating  the two measures in Theorem 3.1 of \cite{zy}, we deduce that
$\Upsilon_1=\Upsilon_2$, $\Gamma_1=\Gamma_2$ which completes the proof of uniqueness.
\vskip 0.3cm

We now prove the inequality \eqref{a0002}.
Define $k:=\|v-\hat{v}\|_\infty^T$ and
$\omega:=(\Xi^{\epsilon_1,\epsilon_2}-\hat{\Xi}^{\epsilon_1,\epsilon_2})-k$.
Then $\omega$ satisfies the following
PDE:
\begin{eqnarray*}
       \frac{\partial{\omega}}{\partial{t}}&=&\frac{\partial^2{\omega}}{\partial{x^2}}  +g^{\epsilon_1,\epsilon_2}(\Xi^{\epsilon_1,\epsilon_2})
       - g^{\epsilon_1,\epsilon_2}( \hat{\Xi}^{\epsilon_1,\epsilon_2}  )
       +\Upsilon^{\epsilon_1,\epsilon_2}-\hat\Upsilon^{\epsilon_1,\epsilon_2}-\big(\Gamma^{\epsilon_1,\epsilon_2}-\hat\Gamma^{\epsilon_1,\epsilon_2}\big)
\end{eqnarray*}
where
 \begin{eqnarray*}
      g^{\epsilon_1,\epsilon_2}(\Xi^{\epsilon_1,\epsilon_2}) & =&\frac{c_1}{[\epsilon_1+ (\Xi^{\epsilon_1,\epsilon_2}+v-\Lambda^1) ]^{\vartheta}}-  \frac{c_2}{[\epsilon_2+ (\Lambda^2-\Xi^{\epsilon_1,\epsilon_2}-v) ]^{\vartheta}}
\end{eqnarray*}
and
\begin{eqnarray*}
      g^{\epsilon_1,\epsilon_2}( \hat{\Xi}^{\epsilon_1,\epsilon_2}  ) & =&\frac{c_1}{[\epsilon_1+ (\hat{\Xi}^{\epsilon_1,\epsilon_2}+\hat{v}-\Lambda^1) ]^{\vartheta}}-  \frac{c_2}{[\epsilon_2+ (\Lambda^2-\hat{\Xi}^{\epsilon_1,\epsilon_2}-\hat{v}) ]^{\vartheta}}.
\end{eqnarray*}

Using Lemma \ref{l1},
we get $\omega_t^+=0$, i.e.
$\Xi^{\epsilon_1,\epsilon_2}-\hat{\Xi}^{\epsilon_1,\epsilon_2}\leq{k}$. By
symmetry,
$\hat{\Xi}^{\epsilon_1,\epsilon_2}-\Xi^{\epsilon_1,\epsilon_2}\leq{k}$. Then
we have for any $T>0$,
\begin{eqnarray*}
\|\Xi^{\epsilon_1,\epsilon_2}-\hat{\Xi}^{\epsilon_1,\epsilon_2}\|_\infty^T\leq\|v-\hat{v}\|_\infty^T.
\end{eqnarray*}
Let $\epsilon_1,\ \epsilon_2\downarrow0$, we conclude
that $\|{\Xi}-\hat \Xi\|_\infty^T\leq\|{v}-\hat v\|_\infty^T$. This is the same as Inequality \eqref{a0002}.
\hfill$\Box$
\vskip 0.3cm
Now we use Picard iteration to get the existence and uniqueness of the solution of Eq. \eqref{a01}.\\
\noindent\textbf{Proof of Theorem \ref{p1}}\quad Similar to Zhang and Yang \cite{zy}, we define
\begin{eqnarray*}
        v_1(x,t)&=&\int_0^1G_t(x,y)X_0(y)dy+\int_0^t\int_0^1G_{t-s}(x,y)f(y,s;X_0)dyds \nonumber\\
        &&+\int_0^t\int_0^1G_{t-s}(x,y)\chi(y,s;X_0)W(dy,ds),
\end{eqnarray*}
where $G_t(x,y)$ is Green's function of the heat equation. So $v_1(x,t)$ satisfies
\begin{eqnarray*}
\left \{\begin{array}{ll}
       \frac{\partial{v_1(x,t)}}{\partial{t}}=\frac{\partial^2{v_1(x,t)}}{\partial{x^2}}+f(x,t;X_0)+\chi(x,t;X_0)\dot{W}(x,t); \\
       v_1(0,t)=v_1(1,t)=0;  \\
      v_1(x,0)=X_0(x) .
\end{array}
\right.
\end{eqnarray*}
Let $(\Xi_1,\Upsilon_1,\Gamma_1)$ be the unique solution of Eq. (\ref{a0001}) with
$v=v_1$. Then
$(X_1,\Upsilon_1,\Gamma_1)$, where $X_1:=\Xi_1+v_1$, is the unique solution of the following
SPDE with two reflecting walls:
\begin{eqnarray*}
\left \{\begin{array}{ll}
       \frac{\partial{X_1(x,t)}}{\partial{t}}=\frac{\partial^2{X_1(x,t)}}{\partial{x^2}}+f(x,t;X_0)+\frac{c_1}{(X_1-\Lambda^1)^{\vartheta}}- \frac{c_2}{(\Lambda^2-X_1)^{\vartheta}}\\
       \ \ \ \ \ \ \ \ \ \ \ \ \ +\chi(x,t;X_0)\dot{W}(x,t)+\Upsilon_1-\Gamma_1; \\
      X_1(0,t)=X_1(1,t)=0;  \\
       X_1(x,0)=X_0(x);\\
       \Lambda^1(x,t)\leq{X}_1(x,t)\leq{\Lambda^2}(x,t).
\end{array}
\right.
\end{eqnarray*}
Iterating this procedure, suppose $X_{n-1}$ has been obtained. Let
\begin{eqnarray*}
      v_n(x,t)&=&\int_0^1G_t(x,y)X_0(y)dy+\int_0^t\int_0^1G_{t-s}(x,y)f(y,s;X_{n-1})dyds \nonumber\\
        &&+\int_0^t\int_0^1G_{t-s}(x,y)\chi(y,s;X_{n-1})W(dy,ds)
\end{eqnarray*}
and $(\Xi_n,\Upsilon_n,\Gamma_n)$ the unique solution of Eq. (\ref{a0001})
replacing $v$ by $v_n$. Then $(X_n,\Upsilon_n,\Gamma_n)$, where $X_n=\Xi_n+v_n$,
is the unique solution of the following SPDE:
\begin{eqnarray*}
\left \{\begin{array}{ll}
       \frac{\partial{X_n(x,t)}}{\partial{t}}=\frac{\partial^2{X_n(x,t)}}{\partial{x^2}}+f(x,t;X_{n-1})+\frac{c_1}{(X_n-\Lambda^1)^{\vartheta}}- \frac{c_2}{(\Lambda^2-X_n)^{\vartheta}}\\
       \ \ \ \ \ \ \ \ \ \ \ \ \
       +\chi(x,t;X_{n-1})\dot{W}(x,t)+\Upsilon_n-\Gamma_n;\\
        X_n(0,t)=X_n(1,t)=0;  \\
      X_n(x,0)=X_0(x);\\
       \Lambda^1(x,t)\leq{X}_n(x,t)\leq{\Lambda^2}(x,t).
\end{array}
\right.
\end{eqnarray*}
From Inequality \eqref{a0002},
\begin{eqnarray}
\|\Xi_n-\Xi_{n-1}\|_\infty^T\leq\|v_n-v_{n-1}\|_\infty^T
\end{eqnarray}
and therefore
\begin{eqnarray}
\|X_n-X_{n-1}\|_\infty^T\leq2\|v_n-v_{n-1}\|_\infty^T.
\end{eqnarray}
Under the assumption (F2,F3), there exists constant $C(k,T)$ depending on $k,T$, such that
\begin{eqnarray*}
\mathbb{E}\big(\|X_n-X_{n-1}\|_\infty^T\big)^k&\leq&2^k\mathbb{E}\big(\|v_n-v_{n-1}\|_\infty^T\big)^k   \nonumber\\
&\leq&{C^{k-1}}(k,T)\mathbb{E}\big(\|X_1-X_0\|_\infty^T\big)^k\frac{T^{n-1}}{(n-1)!}.
\end{eqnarray*}
Hence
there
exists $X(\cdot,\cdot)\in{C}(Q)$
such that
$\underset{n\rightarrow\infty}{lim}\mathbb{E}\big(\|X_n-X\|_\infty^T\big)^k=0$.
Going with the same argument as section 4 in \cite{zy}, we can obtain that
 $(X,\Upsilon,\Gamma)$ is a solution to the SPDE (\ref{a01}). Furthermore, the uniqueness is similar to that in \cite{zy} with the
help of Inequality \eqref{a0002}.
\hfill$\Box$

\section{Pathwise Properties}
\setcounter{equation}{0}
\renewcommand{\theequation}{4.\arabic{equation}}
Let S be the circle [0,1], with the endpoints identified, and let $\rho(x,y)$ be the distance from
$x$ to $y$ along the circle S. That is, let
$$
\rho(x,y)=\min_{k\in\mathbb{Z}}|x-y+k|.
$$
In this section, we
consider pathwise properties of the following SPDEs on $S\times\mathbb{R}_+$:
\begin{equation}\label{4.1}
\left \{\begin{array}{ll}
         \frac{\partial{X}}{\partial{t}}=\frac{\partial^2{X}}{\partial{x^2}}+f(x,t,X)+\frac{c_1}{(X-\Lambda^1)^{\vartheta}}- \frac{c_2}{(\Lambda^2-X)^{\vartheta}}+\chi(x,t,X)\dot{W}+\Upsilon-\Gamma;\\
       X(x,0)=X_0(x)\in{C}(S);\\
       \Lambda^1(x,t)\leq{X}(x,t)\leq{\Lambda^2}(x,t),\ \ {\rm for}\ (x,t)\in{S\times\mathbb{R}_+}.
\end{array}
\right.
\end{equation}\\
Here $\dot{W}(x,t)$ stands for the space-time white noise.

Set $\widehat{\mathcal{A}}_{M,T}=\{X\in C(S\times\mathbb{R}_+):\ \underset{x\in S,t\in[0,T]}{\sup}|X(x,t)|\leq{M}\}$. The reflecting walls $\Lambda^1,\ \Lambda^2$ and  coefficients: $f, \chi$ are assumed to be analogous with that in section 2, i.e.
\\
(H'1)$\Lambda^i(x,t)\in C(S\times \mathbb{R}_+)$ and $\Lambda^1(x,t)<\Lambda^2(x,t)$ for $(x,t)\in S\times\mathbb{R}_+$; \\
(H'2)
$\frac{\partial{\Lambda^i}}{\partial{t}}+\frac{\partial^2{\Lambda^i}}{\partial{x^2}}\in
L^2(S\times[0,T])$;\\
(H'3)
$\frac{\partial{}}{\partial{t}}\Lambda^i(0,t)=\frac{\partial{}}{\partial{t}}\Lambda^i(1,t)=0$
for $t\geq0$;\\
(H'4) $\frac{\partial{}}{\partial{t}}(\Lambda^2-\Lambda^1)\geq0$.\\
and\\
(F'1) $f(\cdot,\cdot;0),\chi(\cdot,\cdot;0)\in L^2(S\times[0,T])$;\\
(F'2) for every $X,\hat{X}\in \widehat{\mathcal{A}}_{M,T}$, there exists $C_{T,M}>0$ such that
$$|f(x,t;X)-f(x,t;\hat{X})|+|\chi(x,t;X)-\chi(x,t;\hat{X})|\leq{C_{T,M}}\underset{y\in S,s\in[0,t]}{\sup}|{X}(y,s)-\hat{X}(y,s)|,$$
for every $x\in S$ and
$t\in[0,T]$;\\
(F'3) there exists $C_T>0$ such that
$$|f(x,t;X)|+|\chi(x,t;X)|\leq{C_T}(1+\underset{y\in S,s\in[0,t]}{\sup}|X(y,s)|),\ \forall (x,t)\in S\times[0,T].$$

Here the precise definition of the solution to Eq. \eqref{4.1} is as follows.
\begin{definition} A triplet
$(X,\Upsilon,\Gamma)$  is a solution to the
Eq.
(\ref{4.1}) if\\
(i) $X=\{X(x,t);(x,t)\in{S\times\mathbb{R}_+}\}$ is a continuous, adapted random field
satisfying $\Lambda^1(x,t)\leq{X}(x,t)\leq{\Lambda^2}(x,t)$,
a.s;\\
(ii) $\Upsilon(dx,dt)$ and $\Gamma(dx,dt)$ are positive and adapted random measures on $S\times\mathbb{R}_+$
satisfying
$$\Upsilon\big(S\times[0,T]\big)<\infty,\ \Gamma\big(S\times[0,T]\big)<\infty$$
for $T>0$;\\
(iii) for all $t\geq0$ and $\psi\in{C^\infty}(S)$ we
have
\begin{eqnarray}
&&\big(X(t),\psi\big)-\big(X_0,\psi\big)-\int_0^t(X(s),\psi^{''})ds-\int_0^t\big(f(y,s,X),\psi\big)ds \nonumber\\
&=&\int_0^t\big(\frac{c_1}{(X-\Lambda^1)^{\vartheta}}
-\frac{c_2}{(\Lambda^2-X)^{\vartheta}},\psi\big)ds+\int_0^t\int_{S}(\chi(y,s,X),\psi)  W(dx,ds) \nonumber\\
&&+\int_0^t\int_{S}\psi\Upsilon(dx,ds)-\int_0^t\int_{S}\psi\Gamma(dx,ds),\ a.s,
\end{eqnarray}
where $(,)$ denotes the inner product in $L^2({S})$;\\
(iv)$$\int_{S\times\mathbb{R}_+}\big(X(x,t)-\Lambda^1(x)\big)\Upsilon(dx,dt)=\int_{S\times\mathbb{R}_+}\big(\Lambda^2(x)-X(x,t)\big)\Gamma(dx,dt)=0.$$
\end{definition}

The existence and uniqueness of the solution of Eq. (\ref{4.1}) is established in a similar way as Theorem \ref{t1}.
\begin{theorem}\label{t4.1}  Let $X_0\in{C}({S})$ satisfy
$\Lambda^1(x,0)\leq{X_0}(x)\leq{\Lambda^2}(x,0),$ for $x\in S$. Under the hypotheses (H'1)-(H'4), (F'1)-(F'3),
there exists a unique solution $(X, \Upsilon,\Gamma)$ to the Eq. \eqref{4.1}
 for all $\vartheta\geq0$.
\end{theorem}

Before studying the pathwise properties of SPDE \eqref{4.1},
we introduce some notation.
In view of assumption (H'2),
suppose that there exist $f_i\in L^2(Q_T)$, where $Q_T:=S\times[0,T]$, such that
$$
\frac{\partial \Lambda^i(x,t)}{\partial t}-\frac{\partial^2\Lambda^i(x,t)}{\partial x^2}=f_i(x,t).
$$
Let
$$
\bar{G}(t,x)=(4\pi t)^{-1/2}\exp(-\frac{x^2}{4t}).
$$
Recall that $\bar{G}(t,x)$ is the fundamental solution of the heat equation on $\mathbb{R}$.
Next, let $G(t,x,y)$ or $G_t(x,y)$ be the fundamental solution of the heat equation on $S$,
and recall that
$$
G(t,x,y)=\sum_{m\in\mathbb{Z}}\bar{G}(t,\rho(x,y)+m).
$$

We have the following
well-known result, the proof is similar to that in Walsh \cite{wa} and is hence omitted.
\begin{lemma}\label{semigroup}
Let $a\in(1,3)$. For any $x,y\in S$ and $t\in[0,T]$, there exists a constant $C_a$ such that
\begin{eqnarray}
\int_0^t\int_{S}|G(t-s,x,y)|^{a}dyds
\leq C_at^{\frac{3a-1}{2}}.
\end{eqnarray}
\end{lemma}

\subsection{The case $\vartheta>3$ }
We consider the special case $\vartheta>3$ in this subsection. We have the following theorem, which is a corollary of Lemma \ref{Lemma 1} below.
The proof of Lemma \ref{Lemma 1} is similar to that in \cite{ml}, but is nontrivial.

\begin{theorem}\label{t2} Let $X_0\in{C}(S)$ satisfy
$\Lambda^1(x,0)<{X_0}(x)<{\Lambda^2}(x,0), x\in S$ and suppose (H'1)-(H'4), (F'1)-(F'3) hold.
If $\vartheta>3$, then
the unique solution $(X, \Upsilon,\Gamma)$ to Eq. (\ref{4.1}) remains in $(\Lambda^1,\Lambda^2)$ for all time and $x\in S$,
hence $\Upsilon=\Gamma=0$.
\end{theorem}
Suppose that
$v(t,x)$ satisfies
\begin{equation}\label{estimate-v-epsilon}
\left \{\begin{array}{ll}
       v_t=v_{xx}+f(x,t,v)+g(x,t,v)+\chi(x,t,v)\dot{W}, \ t>0,\ x\in S;\\
        v(0,x)=X_0(x),
\end{array}
\right.
\end{equation}
where $f,\ \chi$ are given at the beginning of this section and
for every $\delta,\widetilde{\delta}>0$, denote
$$g(x,t,X(x,t))=\Big(|X(x,t)-\Lambda^1(x,t)|\vee(\delta/2)\Big)^{-\vartheta}-\Big(|\Lambda^2(x,t)-X(x,t)|\vee(\widetilde{\delta}/2)\Big)^{-\vartheta}.$$

\begin{lemma}\label{Lemma 1}
Given $T>0$, $\vartheta>3$ and $\epsilon>0$. There exists $\delta_0,\ \delta_1>0$ such that for any
$0<\delta\leq\delta_0$, $0<\widetilde{\delta}\leq\delta_1$, and $X_0(x)\in C(S)$
satisfying $\Lambda^2(x,0)-\widetilde{\delta}\geq X_0(x)\geq \Lambda^1(x,0)+\delta$ for all $x\in S$,
\begin{eqnarray}\label{estimate 1}
\mathbb{P}\Big(\inf_{0\leq t\leq T}\inf_{x\in S}\big( v(x,t)-\Lambda^1(x,t)\big)<\delta/2\Big)<\epsilon,
\end{eqnarray}
and
\begin{eqnarray}\label{estimate 2}
\mathbb{P}\Big(\inf_{0\leq t\leq T}\inf_{x\in S}\big( \Lambda^2(x,t)-v(x,t)\big)<\widetilde{\delta}/2\Big)<\epsilon.
\end{eqnarray}
Here $v(x,t)=\infty$ if $t$ is greater or equal to the blow-up time for $v(x,t)$.
\end{lemma}
\vskip 0.3cm

\noindent\textbf{Proof}\quad
At the beginning of the proof, we introduce some notations. Let $c:=\underset {Q_T} \inf(\Lambda^2(x,t)-\Lambda^1(x,t))$, $C:=\underset {Q_T}\sup(\Lambda^2(x,t)-\Lambda^1(x,t))$.
Denote
\begin{eqnarray*}
\Theta=\{X\in C(Q_T):\ \sup_{(y,s)\in Q_T}|X(y,s)-\Lambda^1(y,s)|\leq c/4\},
\end{eqnarray*}
$$\varphi:=\inf_{\Theta}\inf_{(x,t)\in Q_T}f(x,t,X(x,t))\ \text{and}\ \psi:=\sup_{\Theta}\sup_{(x,t)\in Q_T}f(x,t,X(x,t)).$$

By the assumptions in this section, $c,\ C, \varphi,\ \psi$ are finite constant. For any $l>0$, denote
$$
U_l=\Big\{X\in C(Q_T):\ l/2\leq X(y,s)-\Lambda^1(y,s) \leq 2l,\ \forall (y,s)\in Q_T\Big\}.
$$

Since $\vartheta>3$, it follows that
$
\frac{1}{\vartheta+1}<1/4.
$
Choose $\kappa$ such that $1/(\vartheta+1)<\kappa<1/4$.
Let $\delta_0,\delta_1\leq c/8,$ $0<\delta\leq\delta_0$ and $0<\widetilde{\delta}\leq\delta_1$.

Note that if $X\in U_\delta$, then
\begin{eqnarray*}
\Lambda^2(x,t)-X(x,t)&=&\Lambda^2(x,t)-\Lambda^1(x,t)+\Lambda^1(x,t)-X(x,t)\\
               &\geq&
                     c-2\delta
               \geq
                     c-2\delta_0
               \geq
                     3/4c,
\end{eqnarray*}
and
\begin{eqnarray*}
\Lambda^2(x,t)-X(x,t)&=&\Lambda^2(x,t)-\Lambda^1(x,t)+\Lambda^1(x,t)-X(x,t)\\
               &\leq&
                     \Lambda^2(x,t)-\Lambda^1(x,t)
               \leq
                    C.
\end{eqnarray*}
So we have
$
C\geq |\Lambda^2(x,t)-X(x,t)|\vee(\widetilde{\delta}/2)\geq 3/4 c.
$

Hence, for $a=\varphi-(3/4c)^{-\vartheta}$, $A=\psi-C^{-\vartheta}$,
\begin{eqnarray}\label{A1}
(2\delta)^{-\vartheta}+a\leq f(x,t,X(x,t))+g(x,t,X(x,t))\leq (\delta/2)^{-\vartheta}+A,\ \forall X\in U_\delta.
\end{eqnarray}

Let $w(t,x)$ satisfy
\begin{equation}\label{define w}
\left \{\begin{array}{ll}
w_t(x,t)=w_{xx}(x,t)+f(x,t,w(x,t))+g(x,t,w(x,t))\\
\ \ \ \ \  \ \  \ \ \ \ \ \ +\chi(x,t,w(x,t))\dot{W}(x,t);\\
w(x,k\beta)=\Lambda^1(x,k\beta)+\delta,\ k\beta\leq t< (k+1)\beta,\  x\in S,
\end{array}
\right.
\end{equation}
where
$\beta=2^{-\vartheta-2}\delta^{1+\vartheta}$,
$M\equiv T/\beta$ is an integer and $k=0,\cdots,M-1$.

To get Inequality (\ref{estimate 1}), it is enough to prove
\begin{eqnarray}\label{est-1}
&&\mathbb{P}\Big(
         \inf_{(x,t)\in Q_T}\big(v(x,t)-w(x,t)\big)\geq 0,\
         \inf_{(x,t)\in Q_T}\big(w(x,t)-\Lambda^1(x,t)\big)\geq \delta/2
         \Big)\nonumber \\
         &\geq& 1-\epsilon.
\end{eqnarray}

For $k\beta\leq t<(k+1)\beta$ and
$x\in S$, and assuming that $\Lambda^1(x,s)+\delta/2\leq w(s,x)\leq \Lambda^1(x,s)+2\delta$ for $k\beta\leq s\leq t$ and $x\in S$,
we have
\begin{eqnarray}\label{Eq. w}
w(x,t)=\int_{S}G(t-{k\beta},x,y)\Lambda^1(y,k\beta)dy+\delta+D_k(x,t)+N_{k,L}(x,t),
\end{eqnarray}
where
$L=\underset {X\in U_\delta}\sup \underset {(x,t)\in Q_T} \sup\chi(x,t,X(x,t)),$ and
\begin{eqnarray*}
D_k(x,t)&=&\int_{k\beta}^t\int_{S}G(t-s,x,y)g(y,s,w(y,s))dyds\\
&&+\int_{k\beta}^t\int_{S}G(t-s,x,y)[f(y,s,w(y,s))] dyds,
\end{eqnarray*}
\begin{eqnarray*}
N_{k,L}(x,t)=\int_{k\beta}^t\int_{S}G(t-s,x,y)[\chi(y,s,w(y,s))\wedge L]W(dy,ds).
\end{eqnarray*}
Denote by $\tau_w$ the blow-up time for $|w|$.
Let $w(x,t)=D_k(x,t)=N_{k,L}(x,t)=\infty$ if $t\geq\tau_w$.

Let
$$
\mathcal{N}_k
=\{\omega\in\Omega:
           \Big|\frac{N_{k,L}(x,t)}{(t-k\beta)^\kappa}\Big|\leq \delta^{1-\kappa(1+\vartheta)}2^{-3-2\vartheta+\kappa(\vartheta+2)},
           \forall (t,x)\in(k\beta,(k+1)\beta]\times S\}.
$$

Because of $1-\kappa(1+\vartheta)<0$, if $\delta_0$ is small enough, then from Lemma 2.3 in
\cite{ml}
$\mathbb{P}(\mathcal{N}^c_k)
\leq \frac{\epsilon}{M},\ \forall \delta\in(0,\delta_0],\ k=0,\cdots,M-1.$

Let
$\mathcal{N}=\bigcap_{k=0}^{M-1}\mathcal{N}_k.$ Since
\begin{eqnarray*}
      \mathbb{P}\Big(\inf_{0\leq t\leq T}\inf_{x\in S}\Big[v(x,t)-\Lambda^1(x,t)\Big]<\delta/2\Big)
&\leq&
      \mathbb{P}(\mathcal{N}^c)
\leq
      \sum_{k=0}^{M-1}\mathbb{P}(\mathcal{N}_k^c)
<
    \epsilon,
\end{eqnarray*}
it suffices to show that $w(x,t)\leq v(x,t)$ and $w(x,t)-\Lambda^1(x,t)\geq\delta/2$ for $0\leq t\leq T$
and $x\in S$. By comparison theorem (see Lemma 2.2 in \cite{ml}) and induction, $w(x,t)\leq v(x,t)$ for $0\leq t< k\beta$ and $x\in S$.

Therefore, our aim is to obtain that if the event $\mathcal{N}_k,\ k=0,\cdots,M$ occurs, then,
$
\Lambda^1(x,t)+\delta/2\leq w(x,t)\leq \Lambda^1(x,t)+2\delta,\ k\beta\leq t<(k+1)\beta,\ x\in S, k=0,\cdots,M-1
$ and $\Lambda^1(x,(k+1)\beta)+\delta\leq w((k+1)\beta-,x),\ x\in S$, where
$f(t-)=\underset{s\uparrow t} \lim f(s)$.

Let $t^*$ be the first time
$t\in [k\beta,(k+1)\beta)$ such that for some $x\in S$, $w(x,t^*)=\Lambda^1(x,t^*)+\delta/2$ or $\Lambda^1(x,t^*)+2\delta$. Define $t^*=(k+1)\beta$ if there is no such
time.

Thus, to get inequality (\ref{estimate 1}), we need only verify that on $\mathcal{N}_k$,\\
(i) $t^*=(k+1)\beta$; (ii) $\Lambda^1(x,(k+1)\beta)+\delta\leq w((k+1)\beta-,x),\ x\in S$.

For $x\in S$ and if $t^*<(k+1)\beta$, then by the definition of $\beta$,
\begin{eqnarray*}
      0\leq D_k(x,t^*)
&\leq&
      \beta
          \sup_{X\in U^1_\delta}
          \sup_{(x,t)\in Q_T}\Big(f(x,t,X(x,t))+g(x,t,X(x,t))\Big)\\
&\leq&
      \beta\Big[(\frac{\delta}{2})^{-\vartheta}+A\Big]\\
&=&
     2^{-\vartheta-2}\delta^{1+\vartheta}\Big[(\frac{\delta}{2})^{-\vartheta}+A\Big]=\frac{\delta}{4}+2^{-\vartheta-2}\delta^{1+\vartheta}A
\end{eqnarray*}
and $
    |N_{k,L}(x,t^*)|
\leq
     \delta2^{-3-2\vartheta}.$

Let $I_k(x,t)=\int_{S}G(t-{k\beta},x,y)\Lambda^1(y,k\beta)dyds$. Since $I_k$ and $\Lambda^1$ are the unique solutions of
the following PDEs respectively,
\begin{equation}
\left \{\begin{array}{ll}
\frac{\partial H(x,t)}{\partial t}=\frac{\partial^2H(x,t)}{\partial x^2},\ \ t\in[k\beta, (k+1)\beta),\ x\in S;\\
H(x,k\beta)=\Lambda^1(x,k\beta)
\end{array}
\right.
\end{equation}
and
\begin{equation}
\left \{\begin{array}{ll}
\frac{\partial\widetilde{H}(x,t)}{\partial t}=\frac{\partial^2\widetilde{H}(x,t)}{\partial x^2}+f_1(x,t),\ \ t\in[k\beta, (k+1)\beta),\ x\in S;\\
\widetilde{H}(x,k\beta)=\Lambda^1(x,k\beta),
\end{array}
\right.
\end{equation}
then, by Lemma \ref{semigroup},
\begin{eqnarray}\label{I_k}
  && |I_k(x,t^*)-\Lambda^1(x,t^*)|\nonumber\\
&=&
    |\int_{k\beta}^{t^*}\int_{S}G(t^*-s,x,y)f_1(y,s)dyds|\nonumber\\
&\leq&
    \Big[\int_{k\beta}^{t^*}\int_{S}|G(t^*-s,x,y)|^2dyds\Big]^{1/2}
    \Big[\int_{0}^{T}\int_{S}|f_1(y,s)|^2dyds\Big]^{1/2}\nonumber\\
&\leq&
     \widetilde{C}{(t^*-s)}^{1/4}
\leq \widetilde{C}{\beta}^{1/4}
=
     \widetilde{C}2^{\frac{-\vartheta-2}{4}}\delta^{\frac{1+\vartheta}{4}}.
\end{eqnarray}

Therefore, if $x\in S$ and if $t^*<(k+1)\beta$, then
\begin{eqnarray*}
    |w(x,t^*)-\Lambda^1(x,t^*)-\delta|
&\leq& |I_k(x,t^*)-\Lambda^1(x,k\beta)|+
    |D_k(x,t^*)|+|N_{k,L}(x,t^*)|\\
&<&
     \widetilde{C}2^{\frac{-\vartheta-2}{2}}\delta^{\frac{1+\vartheta}{2}}
    +
     \frac{\delta}{4}+2^{-\vartheta-2}\delta^{1+\vartheta}A
    +
      \delta2^{-3-2\vartheta}.
\end{eqnarray*}

Hence when $\delta_0$ is small enough,
$
|w(x,t^*)-\Lambda^1(x,t^*)-\delta|<\delta/2,\ \forall\delta\in(0,\delta_0],
$ and so $t^*=(k+1)\beta$.

To show (ii), in view of
\begin{eqnarray*}
  w(x,(k+1)\beta-)
&=&
   \int_{S}G((k+1)\beta-k\beta,x,y) \Lambda^1(y,k\beta)  dy  +\delta  \nonumber\\
   &&+D_k(x,(k+1)\beta)+N_{k,L}(x,(k+1)\beta),
\end{eqnarray*}
we have to check
\begin{eqnarray}\label{Eq. N111}
&&|I_k(x,(k+1)\beta)-\Lambda^1(x,(k+1)\beta)|+|N_{k,L}(x,(k+1)\beta)|\nonumber\\
&\leq& D_k(x,(k+1)\beta).
\end{eqnarray}

Therefore, inequality (\ref{Eq. N111}) is acquired by the following inequalities\\
$D_k(x,(k+1)\beta)
\geq
\delta2^{-2\vartheta-2}+a2^{-\vartheta-2}\delta^{1+\vartheta}$,
$|N_{k,L}(x,(k+1)\beta)|
\leq
\delta2^{-2\vartheta-3}$,\\ and
$|I_k(x,(k+1)\beta)-\Lambda^1(x,(k+1)\beta)|
       \leq
           \widetilde{C}2^{\frac{-\vartheta-2}{2}}\delta^{\frac{1+\vartheta}{2}}$.
This completes the proof of (\ref{estimate 1}). Similarly, we can get (\ref{estimate 2}).
\hfill$\Box$

\subsection{The case $0<\vartheta<3$}
In this subsection, we assume the diffusion coefficient $\chi$ is bounded away from zero.
 More precisely, we assume that
$$
\text{(A) There exists $c>0$ such that }c\leq|\chi(x,t,X)|,\ \forall(x,t)\in S\times\mathbb{R}_+,\ X\in\mathbb{R}.
$$
\vskip 0.3cm


Denote $X$ the solution of Eq. (\ref{4.1}), define
$$\tau_1:=  \inf \{ t>0;\ \underset{x\in S}\inf\big( X(x,t)-\Lambda^1(x,t)  \big)=0  \}, $$
$$\tau_2:=  \inf \{ t>0;\ \underset{x\in S}\inf\big( \Lambda^2(x,t)-X(x,t)  \big)=0  \},$$
and set $\tau:=\tau_1\wedge\tau_2$.
%
\vskip 0.3cm

The aim of this subsection is to show that the following proposition holds.
\begin{prop}\label{prop1}
Let the same assumptions as Theorem \ref{t4.1}, and (A) hold. Let $0<\vartheta<3$, and suppose that $X_0\in C(S)$ satisfies $\Lambda^1(x,0)<X_0(x)<\Lambda^2(x,0)$ for all
$x\in S$. If $X$ solves Eq. \eqref{4.1},
then $P(\tau<\infty)>0$, i.e, the solution of Eq. \eqref{4.1} hits $\Lambda^1(x,t)$ or $\Lambda^2(x,t)$ in finite time with positive probability.
\end{prop}

First, we introduce the following SPDE
\begin{equation}\label{b02}
\left \{\begin{array}{ll}
\frac{\partial{v (x,t)}}{\partial{t}}=\frac{\partial^2{v (x,t)}}{\partial{x^2}}+f(x,t,v )+\frac{c_1}{(v -\Lambda^1)^{\vartheta}} +\chi\big(x,t,v \big)\dot{W}(x,t);\\
v (x,0)=X_0(x).
\end{array}
\right.
\end{equation}
Then Eq. \eqref{b02} has a unique solution on the random time interval $[0,\tau_v]$, where
$$\tau_v:=\inf \{ t>0;\ \underset{x\in S}\inf\big( v (x,t)-\Lambda^1(x,t)  \big)\leq0  \}. $$

The next lemma is analogous to that in \cite{mp}.
\begin{lemma}\label{Lem1}
Let $0<\vartheta<3$ and suppose that $X_0\in C(S)$ and $\underset{x\in S}\inf \big(X_0(x)-\Lambda^1(x,0)\big)>0$. If $v$ solves Eq. \eqref{b02},
then $P(\tau_v<\infty)>0$.
\end{lemma}
\noindent\textbf{Proof.}\quad
Let $w(x,t)=v(x,t)-\Lambda^1(x,t)$ and $w(x,0)=X_0(x)-\Lambda^1(x,0)$.
Then $w(x,t)$ is the unique solution of the following equation on $[0,\tau_v]$,
\begin{equation}\label{b03}
\left \{\begin{array}{ll}
\frac{\partial{w (x,t)}}{\partial{t}}=\frac{\partial^2{w (x,t)}}{\partial{x^2}}-f_1(x,t)+f(x,t,w+\Lambda^1 )+\frac{c_1}{w^{\vartheta}}\\
\ \ \ \ \ \ \ \ \ \ \ \ \  +\chi\big(x,t,w+\Lambda^1 \big)\dot{W}(x,t),\ (x,t)\in S\times\mathbb{R}_+;\\
w(x,0)=X_0(x)-\Lambda^1(x,0).
\end{array}
\right.
\end{equation}

Now we use a Girsanov transformation to remove the drift term $-f_1(x,t)+f(x,t,w+\Lambda^1 )$ from Eq. (\ref{b03}).

Consider
\begin{equation}\label{b04}
\left \{\begin{array}{ll}
\frac{\partial{\Xi (x,t)}}{\partial{t}}=\frac{\partial^2{\Xi (x,t)}}{\partial{x^2}}+\frac{c_1}{\Xi^{\vartheta}} +\chi\big(x,t,\Xi+\Lambda^1 \big)\dot{W}(x,t),\ (x,t)\in S\times\mathbb{R}_+;  \\
\Xi(x,0)=X_0(x)-\Lambda^1(x,0).
\end{array}
\right.
\end{equation}
Eq. \eqref{b04} has a unique solution $\Xi$ on $[0,\tau_\Xi]$, where $$\tau_\Xi:=\inf \{ t>0;\ \underset{x\in S}\inf\big( \Xi (x,t)  \big)\leq0  \}.$$

Let $\mathbb{P}_w^{\tau_\Xi\wedge T},\ \mathbb{P}_\Xi^{\tau_\Xi\wedge T}$ denote the measures on the space of continuous function $g:S\times\mathbb{R}_+\rightarrow \mathbb{R}$ induced by $w^{\tau_\Xi\wedge T},\ \Xi^{\tau_\Xi\wedge T}$ respectively. Here
we denote $h^\tau$ the truncated function $h^\tau(x,t)=h(x,t)1_{\{t\leq\tau\}}$. Because of Condition (A), we have $\mathbb{P}_w^{\tau_\Xi\wedge T}$
is absolutely continuous with respect to $\mathbb{P}_\Xi^{\tau_\Xi\wedge T}$, and
$$
\frac{d\mathbb{P}_w^{\tau_\Xi\wedge T}}{d\mathbb{P}_\Xi^{\tau_\Xi\wedge T}}
  =
J
  :=
   \exp\big( \int_0^{T\wedge \tau_\Xi}\int_{S}  q(x,t) W(dx,dt) - 1/2\int_0^{T\wedge \tau_\Xi}\int_{S}  (q(x,t))^2  dxdt \big),
$$
where $q(x,t)=\chi^{-1}(x,t,\Xi+\Lambda^1) [-f_1(x,t)+f(x,t,\Xi+\Lambda^1)] $.

Note the fact that
$$
\int_0^{\tau_\Xi}\int_{S}(q(x,t))^2dxdt<\infty,\ a.s.,
$$
and the fact $\mathbb{P}(\tau_\Xi<\infty)>0$ by Theorem 1 of \cite{mp}, hence $\mathbb{P}(\tau_v<\infty)>0$.
\hfill$\Box$
\vskip 0.3cm
\noindent\textbf{Proof of Proposition \ref{prop1}}\quad
From Theorem 2.1 (comparison theorem) in \cite{dp1}, the solution $X$ of Eq. \eqref{4.1} is less than the solution $v$ of Eq. \eqref{b02} on the time
interval $[0,\tau_v\wedge\tau_1\wedge\tau_2]$, that is
$$
X(x,t)\leq v(x,t),\ \forall x\in S,\ t\in[0,\tau_v\wedge\tau_1\wedge\tau_2].
$$
Hence we have $\{\tau_v <\infty,\  \tau_2 =\infty \}\subset\{\tau_1 <\infty \}$.

By Lemma \ref{Lem1}, we have
\begin{eqnarray*}
0<P(\tau_v <\infty)
     &=&P(\{\tau_v <\infty,\  \tau_2 <\infty \}\cup \{\tau_v <\infty,\  \tau_2 =\infty \} )\\
     &\leq& P(\{\tau_2 <\infty \}\cup \{\tau_1 <\infty \} )\\
     &=& P(\tau_1\wedge\tau_2 <\infty  )\\
     &=& P(\tau <\infty  ).
\end{eqnarray*}
This proves Proposition \ref{prop1}.
\hfill$\Box$

\vskip 0.2cm {\small {\bf  Acknowledgements}\   The authors  thank
Prof. Zhao Dong for his encouragement and
valuable suggestions.

\vskip 0.2cm {\small 

\end{document}